\documentclass[12pt]{article}
\usepackage{amsmath,amssymb,amsthm,amsfonts}
\usepackage{latexsym}
\begin{document}

\newtheorem{prop}{Proposition}[section]
\newtheorem{cor}{Corollary}[section] 
\newtheorem{theo}{Theorem}[section]
\newtheorem{lem}{Lemma}[section]
\newtheorem{rem}{Remark}[section]
\newtheorem{con}{Conjecture}[section]
\newtheorem{as}{Assumption}[section]
\setcounter{page}{1} 
\noindent
\begin{center}
{\Large \bf The  size of the largest
component below phase transition in inhomogeneous random graphs.}
\end{center}

\begin{center}
TATYANA S. TUROVA\footnote{Research was 
supported by the Swedish Natural Science Research
Council.} 
\end{center}

\begin{center}
{\it Mathematical Center, University of
Lund, Box 118, Lund S-221 00, 
Sweden. }
\end{center}

\begin{abstract}
We study the "rank 1 case" of the inhomogeneous random graph model.
In the subcritical case we derive an exact formula for the asymptotic
size of the largest connected component scaled to $\log n$.
This result is new, it completes the corresponding known result in the
supercritical case.
We provide
some
examples of application of a new formula. 
\end{abstract}

\medskip

\medskip

\renewcommand{\theequation}{\thesection.\arabic{equation}}
\section{Introduction.}
\setcounter{equation}{0}

\subsection{Inhomogeneous random graphs.}
We consider here a subclass of a general inhomogeneous
random graph model $G^{\cal V}(n,\kappa )$ with a vertex space
$${\cal V}=(S,\mu, (x_1, \ldots, x_{n})_{n\geq 1})$$
introduced in \cite{BJR}. 
Here $S$ is a separable metric space and $\mu$ is a Borel probability
measure on $S$. 
Recall  the basic definitions
and assumptions from \cite{BJR}. For each $n$ the set of vertices of
the graph $G^{\cal V}(n,\kappa )$ is a deterministic or random
sequence $x_1, \ldots, x_{n}$ of points in $S$, such that for any
$\mu$-continuity set $A\subseteq S$
\begin{equation}\label{set}
    \frac{\#\{i: x_i\in A\}}{n}\stackrel{P}{\rightarrow}   \mu (A).
\end{equation}
Given the sequence $x_1, \ldots, x_{n}$, we let  $G^{\cal V}(n,\kappa)$
 be the random graph on these vertices, such that any
two vertices $x_i$ and $x_j$ are connected by an edge independently of the others
and with a probability
\begin{equation}\label{pe}
p_{x_i,x_j}(n)= \min\{\kappa _{n}(x_i,x_j)/n,1\},
\end{equation}
where $\kappa _{n}$ is a symmetric
nonnegative measurable function on $S\times S$. 
We assume also that for all $x(n)\rightarrow x$ and  $y(n)\rightarrow y$ in $S$
\begin{equation}\label{kap}
    \kappa _{n} (x(n), y(n))\stackrel{a.s.}{\rightarrow}   \kappa (x,
    y)\ \ \
\end{equation}
as $n \rightarrow \infty$, where the kernel
$\kappa$ is {\it graphical} on $\cal V$, which means that

(i) $\kappa$ is continuous $a.s.$ on $S\times S$;

(ii) $ \kappa \in L^1(S \times S,  \mu\times \mu)$;

(iii)
\[\frac{1}{n} {\bf E} e\Big(G^{\cal V}(n,\kappa )\Big)
\rightarrow
\frac{1}{2}\int_{S^2}\kappa (x,y)d {\mu}(x)d
{\mu}(y),\]
where $e(G)$ denotes the number of edges in a graph $G$.

It was observed in \cite{K} that
 random graphs can be naturally related to a certain  branching
process underlying the
algorithm of revealing a connected component in a graph.
This idea was extended  in \cite{T1} for some inhomogeneous graph  model,
where a multi-type branching process was introduced to study the
connectivity of the graph. But it was not until \cite{BJR} that 
a comprehensive theory of
inhomogeneous random graphs was developed, which provided a unified
approach to many models studied previously. 
 
Already in \cite{K} it was shown 
that
in
 the classical random graph model $G_{n,p}$ with $p=c/n$ 
the size of the largest connected component scaled to $n$ 
 asymptotically equals the survival probability of the associated
 branching process.
 A correspondent result was established in \cite{BJR} for a general model
$G^{\cal V}(n,\kappa)$ described above. We shall recall this result here.
Let $ C_1 \Big( G \Big)$ denote the size (the number of vertices) of the largest
connected component in a graph $G$. Then 
Theorem 3.1 from \cite{BJR} states   that
\begin{equation}\label{B}
\frac{C_1 \Big(G^{\cal V}(n,\kappa)
    \Big)}{ n }\stackrel{P}{\rightarrow} \rho_{\kappa}: = \int_{S} \rho_{\kappa}(x)
  d\mu(x),
  \end{equation}
where
$\rho_{\kappa}(x)$
is the survival probability of
a multi-type 
 Galton-Watson process $B_{\kappa}(x)$ defined as
 follows. The type space of  $B_{\kappa}(x)$ is $S$, and initially there
 is a single particle of type $x \in S$. Then 
 at any step,
 a particle of type $x \in S$ is replaced in the next
 generation by a set of particles where the number of  particles of
 type $y$ has a Poisson distribution
with intensity 
 $\kappa
 (x,y)d {\mu}(y)$.
It was also proved in \cite{BJR} that $\rho_{\kappa}(x)$
 is the maximum
  solution to
\[\rho_{\kappa}(x) = 1- e^{-\int_{S} \kappa
 (x,y) \rho_{\kappa}(y)d \mu(y) }. \ \]
 Whether $\rho_{\kappa}$ is zero or strictly positive depends only on
the norm of an integral operator $T_{\kappa}$ defined as 
\begin{equation}\label{Tk}
(T_{\kappa}f)(x)=\int_S{\kappa}(x,y)f(y)d\mu(y)
\end{equation}
with norm
\[\| T_{\kappa}\|= \sup\{\| T_{\kappa}f\|_2:f\geq 0, \|f\|_2\leq1\}.\]
Then according to Theorem 3.1 from \cite{BJR} the survival probability 
\begin{equation}\label{rok}
\rho_{\kappa} \ 
\left\{ \
\begin{array}{ll}
> 0, & \mbox{ if } \ \| T_{\kappa}\|>1, \\
= 0, & \mbox{ if } \ \| T_{\kappa}\|\leq 1.
\end{array}
\right.
\end{equation}
Hence, while (\ref{B}) describes rather accurate the size of the largest
connected component above the phase
transition, i.e., when $\| T_{\kappa}\|>1$, all what we can get from (\ref{B})
 when $\| T_{\kappa}\| \leq  1$ is 
$C_1\Big(G^{\cal V}(n,\kappa )\Big)=o_P(n).$
Only under an additional  assumption
\begin{equation}\label{A2}
\sup_{x,y,n}\kappa _{n}(x,y)<\infty
\end{equation}
Theorem 3.12 in \cite{BJR} proves
in the case $\| T_{\kappa}\| < 1$
 that
 $G^{\cal V}(n,\kappa )=O(\log n)$ {\bf whp} (which means "with high
 probability", i.e., with probability tending to one as $n \rightarrow
 \infty$).

In the case of a homogeneous random graph  $G_{n,p}$ with $p=c/n$
the following
convergence in probability  (and even more precise result)
was derived already in \cite{ER}: if $c<1$ then
\begin{equation}\label{er}
\frac{C_1(G_{n,c/n})}{\log n} \stackrel{P}{\rightarrow}  \frac{1}{c-1+|\log c|}
\end{equation}
as $ n \rightarrow \infty$.
However, the method used in \cite{ER} is not
applicable for an inhomogeneous model.

\subsection{Main results.}

 Our aim here is to
derive the asymptotics
of the size of the largest component scaled to $\log n$ (similar to (\ref{er}))
 for inhomogeneous
random graph model in the case $\| T_{\kappa}\| < 1$.
We show, that this is also directly related to the parameters of the introduced
 branching process $B_{\kappa}$.

Assume from now on that $S \subseteq {\bf R}_+$ is finite or countable, $\mu$ is a
probability on $S$, and 
a graphical kernel $\kappa$ on $S\times S$
has a form
\begin{equation}\label{ka}
\kappa(x,y) =c \psi(x)\psi(y),
\end{equation}
where $\psi $ is a positive function on $S$ and $c$ is a positive constant.
We consider a graph $G^{\cal V}(n,\kappa )$ on the vertex space $\cal
V$ which satisfies condition
(\ref{set}), and  given $x_1, \ldots, x_n,$ the edges are independent and have
probabilities (\ref{pe}) with $\kappa _n=\kappa$, i.e., 
 $$p_{x_i,x_j}(n)= 
\min\left\{\frac{c \psi(x_i)\psi(x_j)}{n},1\right\}.$$
In this case operator $T_{\kappa}$ defined in (\ref{Tk}) has rank 1 (giving the name "the rank 1 case" of inhomogeneous random graph model, see Chapter 16.4 in \cite{BJR}), and 
\begin{equation}\label{F9*}
\| T_{\kappa}\| =c \sum_{S}\psi^2(x)\mu(x).
\end{equation}

\begin{as}\label{Ass1}
Let function  $\psi$
satisfy one of the following conditions: either
\begin{equation}\label{U}
    \sup_{x \in S} \psi(x)<\infty,
\end{equation}
or 
 for some monotone increasing unbounded function $\psi_0$ and
positive constants $A_1 \leq A_2$
\begin{equation}\label{mg}
  A_1\psi_0(x)  \leq \psi (x)\leq A_2\psi_0(x),
\end{equation}
for all large $x$, and
\begin{equation}\label{As1}
\sum_{S}e^{a \psi(x)} \mu(x)< \infty
\end{equation}
for some positive $a$. 
\end{as}

\noindent
Also we shall assume that for any $\varepsilon>0$ and $q>0$
\begin{equation}\label{A3}
{\bf P}\left\{  \left| \frac{\#\{i: x_i=k\}}{n}-\mu(k)\right|\leq
  \varepsilon e^{q\psi(k)}\mu(k), \ \ k \in S \right\} \rightarrow 1
\end{equation}
as $n \rightarrow \infty $. 
Notice that when $S$ is finite, convergence (\ref{A3}) trivially follows by (\ref {set}). 

Let
$B_{\kappa}(x)$ be a branching process defined as above: it starts
with one particle of type $x\in S$, and then
 at any step,
 a particle of type $x \in S$ 
produces $Po\Big(\kappa
 (x,y) {\mu}(y)\Big)$ number of offspring of each type $y\in S$.
 Denote ${\cal X}(x)$ the size of the  total progeny of $B_{\kappa}(x)$, and let
\begin{equation}\label{1}
r(c)=sup \, \{z\geq 1: \sum _S \mu(x) \psi(x){\bf E} z^{ {\cal X}(x)}\  <\infty\}.
\end{equation}

\begin{theo}\label{T1}
Let $\kappa(x,y) = c \psi(x)\psi(y)$ and set
\begin{equation}\label{sc}
 c^{cr}:=\left( \sum_{S}\psi^2(x)\mu(x)\right)^{-1}.
\end{equation}
Under Assumption \ref{Ass1} and (\ref{A3})
we have
\begin{equation}\label{lp}
\frac{C_1 \Big(G^{\cal V}(n,\kappa )
    \Big)}{\log n }
\stackrel{P}{\rightarrow} \frac{1}{\log r(c)}
\end{equation}
as $ n \rightarrow \infty$, where
\begin{equation}\label{r1}
r(c)\left\{
  \begin{array}{ll}
 >1, & \mbox{ if } c<c^{cr}, \\ \\
=1, & \mbox{ if }   c\geq c^{cr}.
\end{array}
\right.
\end{equation}
\end{theo}

Observe that due to (\ref{F9*}) one has
\begin{equation}\label{F9}
\| T_{\kappa}\| < 1 \ \ \ \Leftrightarrow \ \ \ c<c^{cr}.
\end{equation}
Hence, the statement of Theorem \ref{T1} is exactly complementary to (\ref{B}) (under the conditions of Theorem \ref{T1}), since
\[
r(c)>1 \ \mbox{ implies } \ \rho_{\kappa}=0,\]
as well as $ \rho_{\kappa}>0$ implies $  r(c)=1$.
Notice, however, that when $c=c^{cr}$ then both $r(c)=1$ and $\rho_{\kappa}=0 $, and none of statements (\ref{B}) or (\ref{lp}) provides substantial information.

It is rather obvious that one can extend Theorem \ref{T1} for the case of non-countable $S$ under similar assumptions, replacing sum by the integral with respect to $\mu$. It is less apparent, but one
 may conjecture as well, that the statement similar to (\ref{lp}) is not 
restricted to the rank 1 case only.

The rank 1 case is proved to
be versatile for applications. One may interpret $\psi(x)$ as an "activity" of a vertex of type $x$. One particularly often seen choice of $\psi$ is 
 $\psi(x)=x$ on $S=\{1, 2, \ldots\}$. Here "type $x$" can  represent a degree of
 a node as in \cite{BDM} or a size of a macro-vertex as in \cite{TV} (see also Chapter 16.4 in \cite{BJR} on other examples). A special feature of the rank 1 case is that it allows one  to compute $r(c)$ in a rather closed form as we state below. 

\begin{theo}\label{T2}
Assume, the conditions of Theorem \ref{T1} are satisfied.
Let $X$ be a random variable in $S$ with probability function
$\mu$. There exists a unique
$y>1$ which satisfies 
\begin{equation}\label{F12}
y=\frac{1}{c{\bf E} \psi( X)}
\ \frac{{\bf E} \psi( X)\exp \left\{
    c\psi( X)\left( {\bf E} \psi( X) \right)(y-1)\right\} }{
    {\bf E} \psi^2( X)\exp \left\{
    c\psi( X)\left( {\bf E} \psi( X) \right)(y-1)\right\}}.
\end{equation}
Then
\begin{equation}\label{clog}
 r(c)=\frac{1}{c
    {\bf E} \psi^2( X)\exp \left\{
    c\psi( X)\left( {\bf E} \psi( X) \right)(y-1)\right\}} \ .
 \end{equation}
\end{theo} 

\bigskip

Notice that Theorems \ref{T1} and \ref{T2} immediately yield (\ref{er}). Indeed, in the case of a homogeneous model $G_{n, c/n}$ we have $\psi(X)\equiv 1$ in Theorem \ref{T2}, trivially implying $y=1/c$, which together with (\ref{clog})
and (\ref{lp}) gives 
(\ref{er}). 

The result (\ref{lp}) is new for the inhomogeneous random graphs. (In fact, even in the case of $G_{n, c/n}$ the role of 
 $r(c)$ was not  disclosed previously.)
In the subcritical case of  $G_{n, c/n}$  the
method of branching processes was used first in \cite{K} to get
a rough
 upper bound for the  largest connected component. Then in  \cite{T3}
it was shown that for the same case one can get an optimal
upper bound  known from (\ref{er}), again using  the branching processes.

Recent study \cite{TV} indicated that an
upper bound  found there for the largest component in the subcritical case 
of some inhomogeneous random graph model
 should be the optimal one. 
We shall show here that the conjecture from \cite{TV} is a simple corollary of the following modification of Theorem \ref{T1}.

Let $L$ denote a component of a graph $G^{\cal V}(n,\kappa )$. In some applications one studies a function of $L$ in the form 
$\sum_{x_i \in L} \psi(x_i)$, which is natural to call an "activity of a component" if we call $ \psi(x)$ an activity of a vertex of type $x$. A similar
characteristic of a graph, called a "volume"
 was also treated in \cite{CL}.

Consider again the branching processes 
$B_{\kappa}(x)$ defined above, which starts
with one particle of type $x\in S$.
Denote ${\widetilde{\cal X}}(x)$ the set of all the offspring of  the branching processes 
$B_{\kappa}(x)$. (Recall that previously we set $|{\widetilde{\cal X}}(x)|={{\cal X}}(x)$.)
Let $\Phi (x)$ be the sum of "activities" of all the offspring of  the branching processes 
$B_{\kappa}(x)$:
\[\Phi (x) =\sum_{v\in {\widetilde{\cal X}}(x)} \psi(v).\]
This implies that $\Phi (x)$ satisfies the following equality in distribution
\begin{equation}\label{Phi}
\Phi (x) \stackrel{d}{=} \psi(x)+ \sum_{y\in S} \sum_{i=1}^{N_x(y)}\Phi_i(y),
\end{equation}
where $N_x(y) \in $ Po($c\psi(x)\psi(y)\mu(y)$), independent for different $x$ and  $y$; random variables $\Phi (x)$ and $\Phi_i (x), i\geq 1, $ are $i.i.d.$, and also independent for different values of $x$; and a sum over empty set is assumed to be zero.
Define also similar to (\ref{1})
\begin{equation}\label{1*}
\alpha(c)=sup \, \{z\geq 1: \sum _S \psi(x){\bf E} z^{ {\Phi}(x)}\ \mu(x)
 <\infty\}.
\end{equation}
Now we are ready to state another result similar to Theorem \ref{T1}.
\begin{theo}\label{T1*}
Let ${\cal L}$ denote a set of all connected components in  $G^{\cal V}(n,\kappa )$.
Under the conditions of Theorem \ref{T1} 
we have
\begin{equation}\label{lp*}
\frac{\max_{L\in {\cal L}} \sum_{x_i \in L} \psi(x_i)}{\log n }
\stackrel{P}{\rightarrow} \frac{1}{\log \alpha(c)}
\end{equation}
as $ n \rightarrow \infty$, where 
\begin{equation}\label{r1'}
\alpha(c)\left\{
  \begin{array}{ll}
 >1, & \mbox{ if } c<c^{cr}, \\ \\
=1, & \mbox{ if }   c\geq c^{cr}.
\end{array}
\right.
\end{equation}
\end{theo}
One can also find formula for $\alpha(c)$ similar to the one in Theorem \ref{T2}.

\subsection{Example.} 
Let  $G_N(p,c)$ be a graph with the set of vertices $B(N):=\{-N, \ldots, N\}^d$
in $Z^d$, $d\geq 1$, with two types of edges: the short-range edges connect independently 
 with probability $p$ each pair $u$ and $v$  if 
$|u-v|=1$, and the long-range edges
 connect independently any pair of two vertices with
probability $c/|B(N)|$. By this definition there can be none, one or two  edges between two vertices
in
graph $G_N(p,c)$, and in the last case the edges are of different
types. 
Assume, that $0\leq p < p_c(d)$, where $p_c(d)$ is the percolation threshold
in dimension $d$.
As it is shown in \cite{TV}, this graph is  naturally related
to the described above rank 1 case. Consider  the
subgraph of  $G_N(p,c)$ induced by the short-range
edges only, which is a purely bond percolation model.
Let $K_N$ denote the number of open clusters (i.e., connected by
the short-range edges only),
and let 
${\bf X} = \{X_1, X_2, \ldots , X_{K_N}\}$
denote the collection  of all  open clusters $X_i \subseteq
B(N)$. Let also $C$ denote an
open cluster containing the origin.
Recall that 
\begin{equation}\label{P1}
\frac{K_N}{|B(N)|} \ {\rightarrow} \ {\bf
  E}\frac{1}{|C|}
\end{equation}
a.s. and in $L^1$ as $N \rightarrow \infty$ (see, e.g., \cite{G}).
\medskip
Call each set $X_i$ a macro-vertex of type $|X_i|$. Now given a
collection of clusters  ${\bf X}$, introduce another graph $
{\widetilde G}_{N}({\bf X}, p,c)$, whose vertices are macro-vertices
 $X_1, X_2, \ldots , X_{K_N}$.
 The probability that two (macro-)vertices
 $X_i$ 
 and $X_j$ with $|X_i|=x$ and $|X_j|=y$ are connected is derived from the original model
 $G_{N}(p,c)$, which is
\begin{equation}\label{p3}
    {\widetilde p}_{xy}( N) = 1-\left(1-\frac{c}{|B(N)|} \right)^{xy}
    =: \frac{\kappa _{K_N} (x, y)}{K_N}.
\end{equation}
Clearly, the size of the largest connected
component in  $G_{N}(p,c)$ has the following representation
\begin{equation}\label{mC}
C_1\left( G_{N}(p,c)\right)= \max_{L} \sum _{X_i\in L}|X_i|
\end{equation}
where the maximum is taken over all connected components $L$ in $
{\widetilde G}_{N}({\bf X}, p,c)$.

 It was shown in \cite{TV} that graph
 $
{\widetilde G}_{N}({\bf X}, p,c)$ fits the definition of an
inhomogeneous random
graph. In particular, 
\begin{equation}\label{kap*}
    \kappa _{K_N} (x, y)\stackrel{a.s.}{\rightarrow}   \kappa (x,
    y):= c{\bf E}(|C|^{-1})xy\ \ \
\end{equation}
as $N \rightarrow \infty$,  and
\begin{equation}\label{mu} 
  \frac{ \#\{1\leq i \leq K_N: |X_i|=k\}}{K_N} \
\stackrel{a.s.} {\rightarrow} \ \frac{1}{{\bf E}(|C|^{-1})} \  \frac{{\bf P}
  \{|C|=k\}}{k} =:{ \mu}(k)
\end{equation}
 as $N \rightarrow \infty$. (We refer to \cite{TV} for the details.)
It follows also from the results of \cite{TV}, that model 
$
{\widetilde G}_{N}({\bf X}, p,c)$ satisfies
the conditions of Theorem \ref{T1} with $S=\{1, 2,
\ldots \}$, $\kappa(x,y)=c{\bf E}(|C|^{-1})xy$, (we set here $\psi(x)=x$), 
and $\mu$ defined in (\ref{mu}). In this case according to (\ref{Tk})
\[\|T_{\kappa}\|=c{\bf E}(|C|^{-1})\sum_S\psi^2(x)\mu(x)= 
\sum_S x^2\frac{1}{{\bf E}(|C|^{-1})} \  \frac{{\bf P}
  \{|C|=x\}}{x} 
=c{\bf E}|C|, \]
implying that $|T_{\kappa}\|<1$ if and only if $c<\frac{1}{{\bf E}|C|}$. 
Hence, taking into account (\ref{mC}) and also convergence (\ref{P1}) we readily get the result on 
$G_N(p,c)$ model.

\begin{cor}\label{C1}
Assume, that  $d\geq 1$ and  $0\leq p < p_c(d)$.
Then
\[ \frac{C_1\left( G_{N}(p,c)\right)}{\log |B(N)| }
\stackrel{P}{\rightarrow} \frac{1}{\log \gamma(p,c) }
\]
as $ N \rightarrow \infty$, where $\gamma(p,c) = \alpha(c{\bf E}(|C|^{-1}))$
and $\alpha(c)$
is defined by (\ref{1*}) with $\psi(x)=x$, and
\begin{equation}\label{r1*}
\gamma(p,c) \left\{
  \begin{array}{ll}
 >1, & \mbox{ if } c<\frac{1}{{\bf E}|C|}, \\ \\
=1, & \mbox{ if }   c\geq \frac{1}{{\bf E}|C|}.
\end{array}
\right.
\end{equation}
\end{cor}

\hfill$\Box$

\noindent
This result was conjectured in \cite{TV}, where it was proved that
for any $\varepsilon>0$
\[\lim_{N\rightarrow \infty}{\bf P}\left\{\frac{C_1\left( G_{N}(p,c)\right)}{\log |B(N)| } > \frac{1}{\log \gamma(p,c)} +\varepsilon
\right\}=0.\]
We shall also refer to \cite{TV} on more exact description of $\gamma(p,c)$
which is similar to the derivation of (\ref{clog}).

\bigskip
\section{Proofs.}
\subsection{The generating function for the progeny of a branching process.}
Recall that
 ${\cal X}(x)$ denote the total number of the particles
(including the initial one) produced by the 
branching process $B_{\kappa}(x)$, and $\Phi(x)$ is the total activity as defined in (\ref{Phi}). Let for $z\geq 1$
\[h_z(x)={\bf E} z^{{\cal X}(x)}, \ \ \ \ \ 
\ g_z(x)={\bf E} z^{{\Phi}(x)}. \]
Define also 
\[H_z=\sum_{S}\psi(x){\bf E}z^{{\cal X}(x)} \ \mu(x).\]
Then we rewrite (\ref{1})
\[r(c)=\sup \{z\geq 1: H_z <\infty\}.\]
First we shall prove the following lemma, which in particular yields
(\ref{r1}) and (\ref{r1'}).
\begin{lem}\label{L1} Let $\mu$ be a probability on $S$, and let function
  $\psi$ be positive on $S$ and satisfy (\ref{As1}). Write (as in (\ref{sc}))
\[
 c^{cr}:=\left( \sum_{S}\psi^2(x)\mu(x)\right)^{-1}.
\]
Then

\noindent
(I)
\begin{equation}\label{r}
r(c)\left\{
  \begin{array}{ll}
 >1, & \mbox{ if } c<c^{cr}, \\
=1, & \mbox{ if }   c\geq c^{cr};
\end{array}
\right.
\end{equation}
(II)
\begin{equation}\label{rg}
\alpha(c)\left\{
  \begin{array}{ll}
 >1, & \mbox{ if } c<c^{cr}, \\
=1, & \mbox{ if }   c\geq c^{cr};
\end{array}
\right.
\end{equation}
(III) for all $n \geq 1$
\begin{equation}\label{F8}
sup \{z\geq 1: h_z(n) <\infty\}=r(c),
\end{equation}
and
\begin{equation}\label{F8g}
sup \{z\geq 1: g_z(n) <\infty\}=\alpha(c).
\end{equation}
\end{lem}

\noindent
{\bf Proof.}
Note that  
function $h_z(k)$ (as a generating function for a  branching
process)
satisfies the following equation
\[h_z(k) =z \exp {\left\{
\sum_{x\in S} \kappa  (k,x) \mu(x) (h_z(x)-1 )
\right\}}\]
\[=z\exp \left\{
c\psi( k) 
\sum_{x\in S} \psi( x) \mu(x)\left(h_z(x)-1\right)\right\} .\]
Let $X$ denote a random variable in $S$ with distribution $\mu$. Then we can rewrite the last formula as follows
\begin{equation}\label{F15}
\begin{array}{ll}
h_z(k)&  =z \exp \left\{ c   \psi( k)
( H_z-{\bf E} \psi( X))\right\}.
\end{array}
\end{equation}
Multiplying both sides by $\psi( k) \mu (k)$ and summing up over $k$ 
we find for all $z<r(c)$
\begin{equation}\label{F10}
H_z = \sum_{k\in S} \psi( k) \mu (k)z
\exp \left\{ c   \psi( k)
( H_z-{\bf E} \psi( X))\right\}
\end{equation}
\[=z{\bf E} \psi( X)\exp \left\{ c   \psi( X)
( H_z-{\bf E} \psi( X))\right\} .\]
Notice, that
\[H_1= {\bf E} \psi( X),\]
and clearly, $H_z$ is an increasing function of $z$. Hence, equation (\ref{F10}) has a
finite solution for some $z>1$ if and only if the equation
\begin{equation}\label{F11}
  y=\frac{1}{{\bf E} \psi( X)}z{\bf E} \psi( X)\exp \left\{
    c\psi( X)\left( {\bf E} \psi( X) \right) 
(y-1)\right\} .
\end{equation}
(for the same value $z$) has a root  $y>1$. It is easy to see taking into account assumption (\ref{As1}), that at least for some $y>1$ and $z>1$
function on the right in (\ref{F11})
\begin{equation}\label{f}
f(y,z):=zf(y):=z \frac{1}{{\bf E} \psi( X)}{\bf E} \psi( X)\exp \left\{
    c\psi( X)\left( {\bf E} \psi( X) \right) 
(y-1)\right\}
\end{equation}
is
  increasing in both variables, it has all the derivatives of the
second order, and $\frac{\partial ^2}{\partial y^2} f(y,z)>0$. Compute now
\begin{equation}\label{Az2}
\frac{\partial }{\partial  y}f(y,z)|_{y=1, z=1}= 
c {\bf E} \psi^2( X)  =\frac{c}{c^{cr}}.
\end{equation}
Hence, if 
$c\geq c^{cr}$ then for any 
$z>1$
there is no solution $y\geq 1$ to (\ref{F11}).
On the other hand, if $c<c^{cr}$ then there exists $z_0>1$ such that for all
$1\leq z \leq z_0$ there is a finite solution $y\geq 1$ to
(\ref{F11}) which in turn implies existence of $z>1$ for which
(\ref{F10}) has a finite solution $H_z$.
This proves (\ref{r}). The statement (\ref{F8}) follows immediately by
(\ref{F15}).

Exploring formula (\ref{Phi}) it is easy to derive  that  
function $g_z(k)$ 
satisfies the following equation
\[g_z(k) =z^{\psi(k)}\exp \left\{
c\psi( k) 
\sum_{x\in S} \psi( x) \mu(x)\left(g_z(x)-1\right)\right\} .\]
Then with a help of random variable $X$ we can rewrite the last formula as follows
\begin{equation}\label{F15g}
\begin{array}{ll}
g_z(k)&  =z^{\psi(k)} \exp \left\{ c   \psi( k)
( G_z-{\bf E} \psi( X))\right\},
\end{array}
\end{equation}
where
\[G_z=\sum_{S}\psi(x) \mu(x)g_z(x).\]
Multiplying both sides of (\ref{F15g}) by $\psi( k) \mu (k)$ and summing up over $k$ 
we find for all $z<\alpha(c)$
\begin{equation}\label{F10g}
G_z = \sum_{k\in S} \psi( k) \mu (k)z^{\psi(k)}
\exp \left\{ c   \psi( k)
( G_z-{\bf E} \psi( X))\right\}
\end{equation}
\[={\bf E} z^{\psi(X)}\psi( X)\exp \left\{ c   \psi( X)
( G_z-{\bf E} \psi( X))\right\} .\]
The rest of the proof of (\ref{rg}) and (\ref{F8g}) 
is identical to that of (\ref{r}) and (\ref{F8}).
\hfill$\Box$

\subsection{Proof of Theorem \ref{T2}.}
We shall find $r(c)=z_0$ as the (unique!) value for which function $y$ is
tangent to $f(y,z_0)$ (see (\ref{f}) and (\ref{F11})) if $y\geq 1$.

Assume, for some $z_0>1$ the function $y$ is
tangent to $f(y,z)$, and let $y_0$ be the tangency point. Hence, $z_0$
and $y_0$ satisfy the following equations
\[ \begin{array}{rl}
z_0 f'(y_0) & =1,\\ \\
z_0 f(y_0) & =y_0, 
\end{array}
\]
which implies that $y_0$ is the unique solution to 
\[
y_0 = \frac{f(y_0)}{f'(y_0)}=\frac{1}{c{\bf E} \psi( X)}
\ \frac{{\bf E} \psi( X)\exp \left\{
    c\psi( X)\left( {\bf E} \psi( X) \right)(y_0-1)\right\} }{
    {\bf E} \psi^2( X)\exp \left\{
    c\psi( X)\left( {\bf E} \psi( X) \right)(y_0-1)\right\}},
\]
 and then
\begin{equation}\label{F13}
 r(c)=z_0 =\frac{1}{f'(y_0)}=\frac{1}{c
    {\bf E} \psi^2( X)\exp \left\{
    c\psi( X)\left( {\bf E} \psi( X) \right)(y_0-1)\right\}}.
 \end{equation}
 This proves formula (\ref{clog}). \hfill$\Box$

\subsection{Proof of Theorem \ref{T1}.}
We shall assume  here that $\psi$ satisfies (\ref{mg}) and (\ref{As1}).
In the case of (\ref{U}) when $\inf_{x\in S} \psi(x)$ $>0$  one may set $\psi_0(x)=const$ in (\ref{mg}), and the proof will follow by the same argument.
When  (\ref{U}) holds with $\inf_{x\in S} \psi(x)=0$, it is easy to
construct an upper and a lower approximations for the kernel $\kappa$ 
 (consult also \cite{BJR} on approximations), 
so that  the proof will again be reduced to the previous case.

\subsubsection{The lower bound.}
First we shall prove that for any $\delta >0$
\begin{equation}\label{lb}
{\bf P}\left\{  C_1 \Big(
G^{\cal V}(n,\kappa )
    \Big) > \left( \frac{1}{\log r(c)} + \delta \right)\log n \right\}
\rightarrow 0
\end{equation}
as $ N \rightarrow \infty$.

Recall the usual algorithm of finding a
connected component in a random graph. Conditionally on the set of
vertices
${\bf V}:=\{x_1, \ldots , x_n\}$, take any vertex $x_i \in {\bf V}$ to be the root. Find all
the vertices $\{ v^1, v^2,...,v^m\}$ connected to this vertex $x_i$ in the graph $ G^{\cal V}(n,\kappa )$, and
then mark $x_i$ as "saturated". Then for each
non-saturated but already revealed vertex, we find all the vertices
connected to it but which have not been used previously. 
We continue this process until we end
up with a tree of saturated vertices.

Denote  $\tau_{n}(x_i)$ the set of the vertices in
the tree constructed according to the  above algorithm with the
root at a vertex $x_i$.

First we shall prove the following
intermediate result.

\begin{lem}\label{LS} If $c<c^{cr}$ then 
  \begin{equation}\label{S11}
\lim_{n\rightarrow 0}{\bf P} \left\{  C_1\Big( G^{\cal V}(n,\kappa )  \Big) > n^{1/2}
\right\}
=0.
\end{equation}
\end{lem}

\noindent
{\bf Proof.}
Let constant $a$
be the one from the condition (\ref{As1}). Then for any
\begin{equation}\label{q}
0\leq q<a/2
\end{equation}
let us define an auxiliary
probability measure on $S$:
\begin{equation}\label{muq}
  \mu_q(k) = m_q e^{q\psi(k)}\mu(k),
\end{equation}
where the normalizing constant
\[m_q:=\left(\sum_Se^{q\psi(k)}\mu(k) \right)^{-1}>0.\]
Notice that $\mu_0(k) =\mu(k)$ for all $k\in S$, and $m_q$ is continuous on $[0, a/2]$ with 
$m_0=1$.
Fix $\varepsilon>0$ and $0<q<a/2$ arbitrarily and define an event
\begin{equation}\label{calB}
{\cal B}_n =\left\{  \left| \frac{\#\{1\leq i\leq n: x_i=k\}}{n}-\mu(k)\right|\leq
  \varepsilon \mu_q(k), \ \ k \in S \right\}.
\end{equation}

Let $|\tau_{n}(x)|$ denote the
number of vertices in
$\tau_{n}(x)$.
Then we easily derive taking into account the assumption (\ref{A3}) that
\begin{equation}\label{SA18}
 {\bf P} \left\{  C_1\Big(  G^{\cal V}(n,\kappa )   \Big) > n^{1/2}
\right\}
\leq
{\bf P} \left\{ \max _{1\leq i\leq n} |\tau_{n}(x_i)| >  n^{1/2}\mid {\cal B}_n
\right\} +o(1)
\end{equation}
\[
\leq
n\sum_{k\in S}
\left({\mu}(k)+\varepsilon {\mu}_q(k)\right)  
{\bf P} \left\{  
|\tau_{n}(k)| >  n^{1/2}\mid {\cal B}_n 
\right\} +o(1)
\]
as $n \rightarrow \infty$.

To approximate the distribution of $|\tau_{n}(k)|$
we shall use the following branching processes.
Let
$B_{c, q}$ be a process defined similar to
$B_{\kappa}$,
but with
the distribution of the offspring 
\[Po\left(c\psi(x)\psi(y)\mu_q(y)\right)\]
instead of $Po\left(c\psi(x)\psi(y)\mu(y)\right)$. Notice, that
$B_{c, 0}$
is  defined exactly as  $B_{\kappa}$.
We set
\[c^{cr}(q)=\left( \sum_{S}\psi^2(x)\mu_q(x)\right)^{-1}. \]
Clearly, $c^{cr}(q)$ is a continuous function of $q$ on $[0, a/2]$ with
$c^{cr}(0)=c^{cr}$.

Let further ${\cal X}^{c,q}(k)$ denote the total number of the particles
(including the initial one) produced by
the branching process $B_{c, q}$
starting with a single particle
of  type $k$.

\begin{prop}\label{LJ1} For any $c<c^{cr}$ one can find 
$q>0$ and
  $c<c'<\min\{c^{cr}(q), c^{cr}\}$ arbitrarily close to $0$ and $c$,
  correspondingly, such that for all $k\geq 1$ and all large $n$
\begin{equation}\label{SA23}
{\bf P} \left\{ |\tau_n(k) | > n^{1/2} \mid {\cal B}_n
\right\} 
\leq
e^{b_1(log \, n)^4} \, 
{\bf P} \left\{
  {\cal X}^{c',q} (k)> n^{1/2}
\right\} ,
\end{equation}
 where
$b_1 $ is some positive constant independent of $k$ and $n$.
\end{prop}

\noindent
{\bf Proof.}
Observe that at each step of the exploration algorithm which defines $\tau_n$,
the number of the type $y$ offspring 
of a particle of type $x$ has a
binomial distribution $ Bin(N_y',{ p}_{xy}(n))$ where $N_y'$
is the number of the remaining vertices of type $y$.

We shall explore the following 
relation between the binomial and the Poisson  distributions. 
Let $Y_{n,p} \in Bin(n,p)$ and $Z_{\lambda} \in Po(\lambda)$, where $0<p<1/4$ and $\lambda>0$.
Then for all $k\geq 0$
\begin{equation}\label{A21}
{\bf P}
  \{ Y_{n,p}=k\}
\leq (1+\gamma p^2 )^n \, 
{\bf P}\{ Z_{n\frac{p}{1-p}}=k
  \},
\end{equation}
where $\gamma$ is some positive constant (independent of $n$, $k$ and $p$).

Notice that conditionally on ${\cal B}_n$ we have
\begin{equation}\label{Ny}
N_y'\leq \#\{1\leq i\leq n: x_i=y\}\leq n(\mu(y)+\varepsilon \mu_q(y))
\end{equation}
for each $y\in S$. The last inequality implies that for any $y$ such that
\[\#\{1\leq i\leq n: x_i=y\}>0\]
we have
\begin{equation}\label{Fe1}
 n( \mu(y)+\varepsilon \mu_q(y)) \geq 1.
\end{equation}
By the assumptions (\ref{As1}), (\ref{mg}) and (\ref{q}) we have for all large $y$
\[
 \mu(y)+\varepsilon \mu_q(y) \leq
e^{-a\psi(y)}+\varepsilon m_q e^{(q-a)\psi(y)}\leq
 b_2e^{(q-a)A_1\psi_0(y)}\leq
 b_2e^{-aA_1\psi_0(y)/2},
\]
where  $0<b_2<1+2\varepsilon$ for all small $q>0$. 
Combining this with (\ref{Fe1}) we obtain for all large $n$ and $y$ such that
$\#\{1\leq i\leq n: x_i=y\}>0$
\[\frac{1}{n}\leq \mu(y)+\varepsilon \mu_q(y) \leq b_2e^{-aA_1\psi_0(y)/2}.\]
This implies that conditionally on ${\cal B}_n$ 
\[
\max_{x \in \{x_1, \ldots, x_n\}} \psi_0(x)\leq A_3 \log n
\]
for some constant $A_3,$
and thus  conditionally on ${\cal B}_n$ 
\begin{equation}\label{Kn1}
{ p}_{x_i x_j}(n)\leq c \frac{(A_3\log n)^2}{n}
\end{equation}
for all $x_i , x_j \in {\bf V}$.
This and (\ref{Ny}) together with the continuity of $m_q$ allow us 
for any fixed 
positive $\varepsilon _1$
to choose $\varepsilon$ and $q$ in (\ref{calB}) so that
conditionally on ${\cal B}_n$ we get
\[
    N_y' \frac{{ p}_{xy}(n)}{1-{ p}_{xy}(n)}
 \  \leq \left({\mu}(y)+\varepsilon {\mu}_q(y)\right) \frac{n{ p}_{xy}(n)}{1-{ p}_{xy}(n)}
\]

\begin{equation}\label{pB1}
 \leq (1+\varepsilon _1){\mu}_{q}  (y)c\psi (x)\psi(y)=:\mu_{q}  (y)c'\psi (x)\psi(y)
\end{equation}
for all large $n$. In other words, for any
$q>0$ and $c'>c$ arbitrarily close to $0$ and $c$, respectively, and
such that
\begin{equation}\label{A22}
 c<c'<\min\{c^{cr}(q), c^{cr}\},
 \end{equation}
bound (\ref{pB1}) holds for all large $n$.

Now according to (\ref{A21}) and (\ref{pB1}) 
\begin{equation}\label{n2}
{\bf P}\{ 
Y_{N_y' ,{ p}_{xy}(n)} \geq k
\}
\leq 
(
1+\gamma{p}_{xy}(n)^2 
)^{N_y' }\, 
{\bf P}
\{ 
Z_{N_y' 
\frac{{p}_{xy}(n)}{1-{ p}_{xy}(n)}}
\geq k
  \} 
\end{equation}
\[
\leq (1+\gamma{ p}_{xy}(n)^2 )^{n}\, 
{\bf P} \{ Z_{\mu_{q}  (y)c'\psi (x)\psi(y)}\geq k
  \}.
\]
Hence, if conditionally on ${\cal B}_n $
 at each (of at most $n$) step
of the exploration algorithm which reveals $\tau_n(k)$, we replace the
$ Bin(N_y',{ p}_{xy}(n))$ variable 
with the
$Po\left(\mu_{q}  (y)c' \psi (x)\psi(y)\right)$ one, we arrive at
the
following bound using branching process $B_{c',q} (k)$ and bound (\ref{Kn1}):
\begin{equation}\label{S12}
{\bf P}
\left\{ |\tau_n(k) | > n^{1/2}\mid {\cal B}_n
\right\}
\leq
\left(1+\gamma \left(c\frac{(A_3\log n)^2}{n}
\right)^2\right)^{n^2}
 \, {\bf P} \left\{ {\cal X}^{c',q} (k)> n^{1/2}
\right\} .
\end{equation}
This implies statement (\ref{SA23}) of the
 Proposition. \hfill$\Box$

\bigskip

Substituting (\ref{SA23})  into (\ref{SA18})
we derive for any
$q>0$ and $c'>c$ that
\begin{equation}\label{SA24}
{\bf P} \left\{   C_1\Big(  G^{\cal V}(n,\kappa )  \Big) >n^{1/2}\right\} 
\leq
b_3 n e^{b_1(\log n)^4}
\sum_{k\in S} 
 \mu _{q}(k)
{\bf P} \left\{ {\cal X}^{c',q} (k)> n^{1/2}
\right\}
 +o(1)
\end{equation}
as $n \rightarrow \infty$, where
$b_3$ is some positive constant. By the Markov's inequality
\begin{equation}\label{SMark}
{\bf P} \big\{{\cal X}^{c',q} (k) > n^{1/2}\big\} \leq z^{-n^{1/2}}
{\bf
    E} z^{{\cal X}^{c',q} (k)}
\end{equation}
for all $z\geq 1$. This bound together with  (\ref{SA24}) yield
\begin{equation}\label{SA25}
{\bf P} \left\{   C_1\Big(  G^{\cal V}(n,\kappa )  \Big) >n^{1/2}\right\} 
\leq b_3 n e^{b_1(\log n)^4}z^{-n^{1/2}}
\sum_{k\in S} 
 \mu _{q}(k){\bf
    E} z^{{\cal X}^{c',q} (k)}
+o(1).
\end{equation}
We are left to show that for some $z>1$ 
\begin{equation}\label{Ma1}
\sum_{k\in S} 
 \mu _{q}(k){\bf
    E} z^{{\cal X}^{c',q} (k)}<\infty .
  \end{equation}
  Note that $\mu_q$ and $\psi$ satisfy the conditions of Lemma
  \ref{L1} for any $0\leq q <a/2$, and moreover $c'<c^{cr}(q)$. Hence, by Lemma \ref{L1} there exists
  $z_0>1$ such that
  \[\sum_{k\in S} 
 \psi(k)\mu _{q}(k){\bf
    E} z_0^{{\cal X}^{c',q} (k)}<\infty,\]
  which clearly implies (\ref{Ma1}), and the statement of Lemma
  \ref{LS} follows by (\ref{SA25}) where we set $z=z_0>1$. \hfill$\Box$

  \bigskip

Now we are ready to complete the proof of (\ref{lb}), following almost the same arguments as in the proof of the previous lemma.
Denote
\[ {\cal B}'_n := {\cal B}_n \cap
\left( C_1\Big( G^{\cal V}(n,\kappa ) \Big)
\leq n^{1/2}
\right).
\]
According to assumption (\ref{A3}) and Lemma \ref{LS}
we have
\[{\bf P} \left\{ {\cal B}'_n\right\}=1- o(1)\]
as $n \rightarrow  \infty$.
This allows us to  derive similar to (\ref{SA18}) for any $\omega$
\begin{equation}\label{A18}
 {\bf P} \left\{  C_1\Big(  G^{\cal V}(n,\kappa )   \Big) > \omega
\right\}
\leq
{\bf P} \left\{ \max _{1\leq i\leq n} |\tau_{n}(x_i)| > \omega \mid {\cal B}_n'
\right\} +o(1)
\end{equation}
\[
\leq
n\sum_{k\in S}
\left({\mu}(k)+\varepsilon {\mu}_q(k)\right)  
{\bf P} \left\{  
|\tau_{n}(k)| > \omega \mid {\cal B}'_n 
\right\} +o(1)
\]
as $n \rightarrow \infty$.
Repeating the same argument which led to (\ref{S12}),
we get the
following bound using the introduced branching process:
\[
{\bf P} \left\{|\tau_{n}(k)| > \omega  \mid {\cal B}'_n
\right\} 
\leq
\left(1+\gamma \left(c\frac{(A_3 \log n)^2}{n}
\right)^2\right)^{b_1 n \sqrt{n}}
 \, {\bf P} \left\{ {\cal X}^{c',q} (k)> \omega
\right\} 
\]
as $n \rightarrow \infty$, where we take into account that we can
perform at most $\sqrt{n}$ steps of exploration (the maximal possible
number of  macro-vertices in any connected component
conditioned on ${\cal B}'_n$ ). Notice also that we can choose here
$c'$ and $q$ arbitrarily close to $c$ and $0$, correspondingly, and so
that condition (\ref{A22}) is fulfilled. 
The last bound  implies
\begin{equation}\label{A23}
{\bf P} \left\{ |\tau_{n}(k)| > \omega \mid {\cal B}'_n 
\right\} 
\leq
(1+o(1))
{\bf P} \left\{  {\cal X}^{c',q} (k)> \omega
\right\} 
\end{equation}
as $n\rightarrow \infty$. Substituting (\ref{A23}) into (\ref{A18})
we derive
\begin{equation}\label{A24}
   {\bf P} \left\{  C_1\Big(  G^{\cal V}(n,\kappa )   \Big) > \omega
\right\}
\leq
 b
n
\sum_{k\in S} 
  \mu _{q}(k)
{\bf P} \left\{ {\cal X}^{c',q}(k)> \omega 
\right\}
 +o(1)
\end{equation}
as $n\rightarrow \infty$, where $b$ is some positive constant. 
Then 
similar to (\ref{SA25})
 we derive from (\ref{A24}) 
\begin{equation}\label{A25}
{\bf P} \left\{  C_1 \Big(G^{\cal V}(n,\kappa ) \Big) > \omega \right\} 
\leq
b
n z^{-\omega} \sum_{k\in S} 
 \mu _{q}(k)
{\bf E} z^{{\cal X}^{c',q} (k)} + o(1).
\end{equation}
Since $\psi$ and probability $\mu_q$
satisfy the conditions of Lemma \ref{L1}, for all values $c'$ and $q$ for which (\ref{A22}) holds, we have by  this lemma
\[r_q(c'):=\sup \{z\geq 1: \sum_{x\in S}\psi(x){\bf E}z^{{\cal X}^{c',q} (x)} \ \mu_q(x) <\infty\} >1,\]
and therefore 
for all $1<z < r_q(c')$ 
\begin{equation}\label{Ma2}
\sum_{k\in S} 
 \mu _{q}(k)
 {\bf E} z^{{\cal X}^{c',q} (k)} <\infty .
\end{equation}
Note that $r_q(c')$ is a continuous function of $q$ and $c'$: this can be explicitly derived from 
formula (\ref{F13}) where $X$
is replaced by a random variable $X_q$ with a probability function $\mu _q$.
Moreover, $r_0(c)=r(c)$. Hence, for any $\delta >0$ 
we can choose a small $\delta' >0$ and $(q,c')$
close to $(0,c)$ so that (\ref{Ma2}) holds with 
\begin{equation}\label{Ma6}
z=r_q(c')-\delta' > 1,
\end{equation}
and also
\begin{equation}\label{Ma4}
\left(\frac{1}{\log r(c)}+\delta \right)\log (r_q(c')-\delta') > 1.
\end{equation}
Now setting
$\omega= \left(\frac{1}{\log r(c)}+\delta \right)\log n$
and $z=r_q(c')-\delta'$
in (\ref{A25}) we derive with a help of (\ref{Ma2})
\[
{\bf P} \left\{  C_1 \Big(G^{\cal V}(n,\kappa ) \Big) > \left(\frac{1}{\log r(c)}+\delta \right)\log n \right\} 
\leq b_3nz^{-\omega}+o(1)
\]
\begin{equation}\label{Ma3}
=
b_3 n \exp\{-\log (r_q(c')-\delta')\left(\frac{1}{\log r(c)}+\delta \right)\log n\}
 + o(1) 
\end{equation}
where $b_3$ is some finite positive constant. This together with (\ref{Ma4})
clearly implies  statement (\ref{lb}). \hfill$\Box$

\subsubsection{The upper bound.}
Here we show 
that for any $\delta >0$
\begin{equation}\label{ub}
\lim _{n \rightarrow \infty}{\bf P}\left\{  C_1 \Big(
G^{\cal V}(n,\kappa )
    \Big) < \left( \frac{1}{\log r(c)} - \delta \right)\log n \right\}
= 0.
\end{equation}

Given a random graph $G^{\cal V}(n,\kappa )$
let $T$ be the number of its connected components, and denote
$W_i$, $i=1, \ldots, T$, the corresponding sets of the vertices in
these components, ordered
arbitrarily.
Let us fix $\varepsilon>0$ arbitrarily and
for each $n \geq 1$  introduce an event
\[{\cal A}(n)= \left\{
\max_{1\leq i\leq T}
|W_i|\leq \left( \frac{1}{\log r(c)} + \varepsilon \right)
\log n \right\} \cap {\cal B}_n\]
with ${\cal B}_n$ defined in (\ref{calB}).
According to (\ref{lb}) and the assumption (\ref{A3})
\begin{equation}\label{F2}
{\bf P} \left\{{\cal A}(n)\right\}\rightarrow 1
\end{equation}
as $ n \rightarrow \infty$.
Let
\[{\bf P}_{{\cal A}(n)}(\cdot)= {\bf P} \left\{\cdot \mid {\cal A}(n) \right\}\]
denote the conditional probability. 
Denote further
\begin{equation}\label{om}
\omega= \left(\frac{1}{r(c)}-\delta\right)\log n, \ \ \
\omega_1= \left(\frac{1}{r(c)}+\varepsilon\right)\log n,
\end{equation}
\begin{equation}\label{fN}
N=N(n)= \frac{n}{\omega_1^2}.
\end{equation}
Clearly, 
\begin{equation}\label{F3}
{\bf P}_{{\cal A}(n)} \left\{T \geq [N]+1\right\}= 1.
\end{equation}

We shall reveal recursively  $[N]+1$ connected components in the graph $G^{\cal V}(n,\kappa )$ in the following way.
Let $V_1$ be a random vertex uniformly 
distributed on ${\bf V}$. Set $L_1=\tau_n(V_1)$ to be
the set of the vertices in the connected component containing vertex $V_1$.

Further for any $U \subset {\bf V} $ let $\tau_n^{U}(v)$
denote a set of vertices of the tree constructed in the same way as 
$\tau_n(v)$ but on the set of vertices ${\bf V}\setminus U$ instead of
${\bf V}$. In particular, with this notation
$\tau_n^{\emptyset}(v)=\tau_n(v)$.

Given  constructed components $L_1, \ldots , L_k$ for  $1 \leq k \leq
[N]$, let 
$V_{k+1}$ be a vertex uniformly distributed on
$V\setminus \cup_{i=1}^kL_i$, and set 
$L_{k+1}=\tau_n^{\cup_{i=1}^kL_i }(V_{k+1})$.
Then according to (\ref{F3}) and (\ref{F2})  we have
\begin{equation}\label{F1}
{\bf P}\left\{  C_1 \Big(
G^{\cal V}(n,\kappa )
    \Big) < \left( \frac{1}{\log r(c)} - \delta \right)\log n \right\}\leq
{\bf P}_{{\cal A}(n)} \left\{\max_{1\leq i\leq [N]+1}|L_i|<\omega \right\}
+o(1)
\end{equation}
as $n \rightarrow \infty$.

Let $x_0\in S$ be such that $\psi(x_0)= \min_{x\in S}\psi(x)$. Then a vertex of type $x_0$ has
 among all different types $x\in S$ the smallest
probabilities of the incident edges, which are $c\psi(x_0)\psi(y)/n$,
$y\in S$. This  implies
that
for any $U\subset {\bf V}$ the size of $\tau_n^U(x_0)$ is stochastically
dominated by $|\tau_n^U(x)|$ for any $x\in S$. 
Notice also that if $U \subset U' $ then  $|\tau_n^{U'}(x)|$ is
stochastically
dominated by $|\tau_n^U(x)|$ for any $x\in S$. This allows us to derive the 
following bound
\begin{equation}\label{F4}
{\bf P}_{{\cal A}(n)} \left\{\max_{1\leq i\leq [N]+1}|L_i|<\omega \right\}
\leq \left(\max_{U\subset {\bf V}: |U|\leq  N\omega_1 }
{\bf P}_{{\cal A}(n)} \left\{|\tau_n^U(x_0)|<\omega \right\}
\right)^{N}.
\end{equation}

To approximate the distribution of $|\tau_n^{U}(x)|$ we
introduce another branching process which will be stochastically
dominated by $B_{\kappa}$. 
First define for any value $D\in S$ another auxiliary probability measure $\hat{\mu}_D$
\begin{equation}\label{muD}
\hat{\mu}_D(y)=
\left\{
\begin{array}{ll}
M_D^{-1}\mu(y), & \mbox{ if } y\leq D, \\
0, & \mbox{ otherwise },
\end{array}
\right.
\end{equation}
where $M_D:=\sum_{y\leq D}\mu(y)$ is a normalizing constant. 
Then for any positive $c$
and $D$ let $\hat{B}_{c, D}$ be a process defined similar to
$B_{\kappa}$,
but with
the distribution of the offspring 
\begin{equation}\label{BD}
Po\left(c\psi(x)\psi(y)\hat{\mu}_D(y)\right)
\end{equation}
instead of $Po\left(c\psi(x)\psi(y)\mu(y)\right)$. Notice, that
$\hat{B}_{c, \infty}$
is  defined exactly as  $B_{\kappa}$.
Let $\hat{\cal X}^{c,D}(x)$ denote the total number of the particles
(including the initial one) produced by
the branching process $\hat{B}_{c, D}$ 
starting with a single particle
of  type $x$.

\begin{lem}\label{L2}
For any $c'<c$ there exists  finite $D$  such that
\begin{equation}\label{F7}
{\bf P}_{{\cal A}(n)}
 \left\{|\tau_n^{U}(x_0)|<\omega \right\}\leq 
\left(1+b\frac{\log^4n}{n^2}\right)^{n\omega_1}
 \, {\bf P} \left\{ \hat{\cal X}^{c',D} (x_0)<\omega 
\right\} 
\end{equation}
for all large $n$ uniformly in $U \subset {\bf V}$ with $|U|\leq N\omega_1
=\left(\frac{1}{r(c)}+\varepsilon\right)^{-1}n/\log n$, where $b=b(c)$ is some positive constant
independent of $c'$
and $D$.
\end{lem}

\noindent{\bf Proof.}
At each step of the exploration algorithm which defines $\tau_n^{U}(x_0)$,
the number of the type $y$ offspring 
of a particle of type $x$ has a
binomial distribution $ Bin(N_y',{ p}_{xy}(n))$ where $N_y'$
is the number of remaining vertices of type $y$. 
We shall first find a lower bound for $N_y'$.

According to assumption (\ref{set}) for any $D\in S$ and $\varepsilon_1>0$
there exists $n(D, \varepsilon_1)$
such that 
\begin{equation}\label{Ma12}
  N_y:=\#\{x_i \in {\bf V}: x_i=y\} \geq (1-\varepsilon_1)\mu(y)n
\end{equation}
for all $y<D$ and $n\geq n(D, \varepsilon_1)$.
Note that conditionally on ${\cal A}(n)$ the number of vertices in
$\tau_n^{U}(x)$ is at most $\omega_1$,
and by deleting an  arbitrary set  $U$ with $|U|\leq N\omega_1$
from ${\bf V}$, we may delete at most $N\omega_1$ vertices of type $y$.
Hence,  at any step of the 
exploration algorithm which defines $\tau_n^{U}(x)$,
 the number $N_y'$ of the remaining vertices of type $y$, is bounded
from below as follows
\[
N_y' \geq N_y-\omega_1 -N\omega_1,
\]
and thus
according to (\ref{Ma12}) 
\[
N_y' \geq n(1-\varepsilon_1)\mu(y)-\omega_1 -N\omega_1
\]
for all $y<D$ and $n\geq n(D, \varepsilon_1)$.
Taking into account definitions (\ref{fN}) and (\ref{om}) we derive from here that for any $\varepsilon _1>0$ and $D>0$
there exists $n(D, \varepsilon_1)$
such that 
\[
  N_y' \geq (1-\varepsilon_1)\mu(y)n
\]
for all $y\leq D$ and $n\geq n(D, \varepsilon_1)$. This implies that
conditionally on ${\cal A}(n)$  at any step of the exploration
algorithm we have
\begin{equation}\label{Ma13}
  N_y' 
\frac{{p}_{xy}(n)}{1-{ p}_{xy}(n)} \geq
\mu  (y)(1- \varepsilon_1)c\psi (x)\psi(y)
\end{equation}
for any fixed $\varepsilon _1, D$, and all $n \geq n(D,\varepsilon_1)$ and $y \leq D$. 
Now with a help of (\ref{muD}) we rewrite (\ref{Ma13})
as follows: 
\begin{equation}\label{Ma14}
  N_y' 
\frac{{p}_{xy}(n)}{1-{ p}_{xy}(n)} \geq \hat{\mu}_D(y)M_D(1- \varepsilon_1)c\psi (x)\psi(y) =: \hat{\mu}_D(y)c'\psi (x)\psi(y) 
\end{equation}
for any fixed $\varepsilon _1, D$, and all $n \geq n(D,\varepsilon)$ and $x,y \in S$, where
\[c'=M_D(1- \varepsilon_1)c.\]
Recall that $\lim_{D\rightarrow \infty}M_D \uparrow 1$. Therefore choosing appropriately constants $D$ and $\varepsilon_1$ we can make $c'$ arbitrarily close to $c$.
Using again  
relation (\ref{A21}) between the Poisson and the binomial distributions, and taking into account  (\ref{Ma14}), we derive for all $k$
\begin{equation}\label{Ma9}
{\bf P}\{ 
Y_{N_y' ,{ p}_{xy}(n)} \leq k
\}
\leq (1+ \gamma{ p}_{xy}^2(n) )^{N_y' } \, 
\, 
{\bf P}
\{ 
Z_{N_y' 
\frac{{p}_{xy}(n)}{1-{ p}_{xy}(n)}}
\leq k
  \} 
\end{equation}
\[
\leq (1+\gamma{ p}_{xy}^2(n) )^{n}\, 
{\bf P} \{ Z_{\hat{\mu}_D(y)c'\psi (x)\psi(y)}\leq k
  \}.
\]
This implies that if conditionally on ${\cal A}_n $,
 at each of at most $\omega_1$ steps
of the exploration algorithm which reveals $\tau_n^{U}(x_0)$, we replace the
$ Bin(N_y',{ p}_{xy}(n))$ variable 
with the
$$Po\left(\hat{\mu}_D(y)c' \psi (x)\psi(y)\right)$$ one, we arrive at
the
following bound using the branching process $\hat{B}_{c, D}$ and
bound (\ref{Kn1})):
\begin{equation}\label{Ma15}
{\bf P}_{{\cal A}(n)}
 \left\{|\tau_n^{U}(x_0)|<\omega \right\}
\leq
\left(1+\gamma c^2A_3^4\frac{\log^4n}{n^2}
\right)^{n\omega_1}
 \, {\bf P} \left\{ \hat{\cal X}^{c',D} (x_0)<\omega 
\right\} 
\end{equation}
which holds for all large $n.$ 
This implies the statement of Lemma \ref{L2}.
\hfill$\Box$

Lemma \ref{L2} together with (\ref{F4}) implies that for all large $n$
\[
{\bf P}_{{\cal A}(n)} \left\{\max_{1\leq i\leq [N]+1}|L_i|<\omega \right\}
\leq \left(\left(1+b\frac{\log^4n}{n^2}
\right)^{n\omega_1}
{\bf P} \left\{ \hat{\cal X}^{c',D} (x_0)<\omega 
\right\} 
\right)^{N}.\]
Substituting this into (\ref{F1})
we derive
\[
{\bf P}\left\{  C_1 \Big(
G^{\cal V}(n,\kappa )
    \Big) < \omega \right\}\leq
\left(\left(1+b\frac{\log^4n}{n^2}
\right)^{n\omega_1}
{\bf P} \left\{ \hat{\cal X}^{c',D} (x_0)<\omega 
\right\} 
\right)^{N}
+o(1)\]
\[=\left(1+b\frac{\log^4n}{n^2}
\right)^{n^2/\omega_1}
{\bf P} \left\{ \hat{\cal X}^{c',D} (x_0)<\omega 
\right\}^{n/\omega_1^2}
+o(1)\]
\begin{equation}\label{F18}
\leq e^{ b_1\log^3n}
\, \Big(1-{\bf P} \left\{ \hat{\cal X}^{c',D} (x_0)\geq\omega 
\right\} 
\Big)^{n/\omega_1^2} +o(1)
\end{equation}
as $n \rightarrow \infty$, where $b_1=b_1(c)$ is a positive constant which depends only on $c$.

Since $\hat{\mu}_D$ and $\psi$ obviously satisfy the conditions of Lemma
\ref{L1}
we get by (\ref{F8}) and (\ref{r}) that if
\[c'<{\hat c}^{cr}(D):=\left(\sum_S \psi^2(x){\hat \mu}_D(x)\right)^{-1},\]
then
\[{\hat r}(c', D):=\sup\{z\geq 1: {\bf E}
z^{\hat{\cal X}^{c',D} (x_0)} <\infty \} >1,\]
and ${\hat r}(c', D)$ can be derived from  (\ref{clog}) where 
$X$ is replaced by a random variable on $S$ with a distribution 
${\hat \mu}_D.$ It is clear that
\[\lim_{D \rightarrow \infty}{\hat c}^{cr}(D)={c}^{cr},\]
${\hat r}(c, D)$ is continuous in $c$, and 
\[\lim_{D \rightarrow \infty}{\hat r}(c, D)=r(c).\]
Hence we can find 
for any given $\delta_1>0$ a large constant $D$ and
$c'<c$ sufficiently close to $c$, such that 
\[{\hat r}(c', D)<r(c)+\delta_1/2.\]
Then it follows from the definition of ${\hat r}(c', D)$ that for some 
positive constant $A=A(\delta_1)<\infty$ and any positive $\omega$ 
\begin{equation}\label{F16}
{\bf P}
 \left\{\hat{\cal X}^{c',D}(x_0)>\omega \right\} \geq 
A({\hat r}(c', D)+\delta_1/2 )^{-\omega } \geq 
A({r}(c)+\delta_1 )^{-\omega }.
\end{equation}
This allows us to  derive from (\ref{F18}) that for any $\delta >0$,
$\delta_1>0$ and some  positive $A$
\begin{equation}\label{F17}
{\bf P}\left\{  C_1 \Big(
G^{\cal V}(n,\kappa )
    \Big) < \omega \right\}
\end{equation}
\[\leq e^{ b_1\log^3n}
 \left(1-A(r(c)+\delta_1 )^{-\left(\frac{1}{\log
        r(c)}-\delta\right) \log n }\right)^{\alpha^2 n/\log^2 n } +o(1)
\]
where $\alpha= \left(\frac{1}{\log
        r(c)}+\varepsilon\right)^{-1}$.
Now for  any $\delta >0$ we choose a positive $\delta_1$ so that
\[\gamma_1:= \left(\frac{1}{\log
        r(c)}-\delta\right) \log \left(r(c)+\delta_1\right) <1.\]
Then (\ref{F17}) becomes
    \begin{equation}\label{F19}
{\bf P}\left\{  C_1 \Big(
G^{\cal V}(n,\kappa )
    \Big)<\left(\frac{1}{\log r(c)}-\delta\right) \log n\right\}
\leq 
e^{b_1\log^3n} \left(1-\frac{A}{ n^{\gamma_1} }\right)^{
    \alpha^2 n/\log^2 n } +o(1),
\end{equation}
where the right-hand side goes to zero when $n\rightarrow
\infty$. This
completes the proof of (\ref{ub}), which together with  (\ref{lb})
yields the assertion of Theorem \ref{T1}. \hfill$\Box$

\subsection{Proof of Theorem \ref{T1*}.}

The proof of Theorem \ref{T1*} almost
exactly repeats the proof of 
Theorem \ref{T1}. The only difference is that a random
variable $|\tau_n(x)|$ used in the proof of 
Theorem \ref{T1} should be replaced by
\[\Psi_n(x):= \sum_{v\in \tau_n(x)}\psi(v),\]
while ${\cal X}^{c,q}$ and $\hat{\cal X}^{c,D}$ should be replaced by
$\Phi^{c,q}$
and $\hat{\Phi}^{c,D}$ which denote the activity of the total progeny of the branching processes $B_{c,q}$ and $\hat{B}^{c,D}$, correspondingly
(see definition of $\Phi$ and (\ref{Phi})). Then due to the results (\ref{F8g})
and (\ref{rg}) on $\alpha$
from Lemma \ref{L1} the proof of Theorem \ref{T1*} follows exactly the
same lines as the proof of Theorem \ref{T1}. \hfill$\Box$

\bigskip

\bigskip

\noindent
{\bf Acknowledgment} The author thanks A. Martin-L{\"o}f for the
helpful discussions.

\end{document}